%% file: main.tex
\documentclass[11pt]{article}
\usepackage{algorithm,algorithmic,amsmath,amssymb,epsfig}
\usepackage[a4paper]{geometry}
\usepackage{url}
\usepackage{layout}

\usepackage{color}
\usepackage{gauss}
\usepackage{graphicx}
\usepackage{epsfig}
\usepackage{subfig}
\usepackage{tikz}

\newcommand{\bfX}{{\textbf{X}}}
\newcommand{\bfY}{{\textbf{Y}}}
\newcommand{\bfB}{{\textbf{B}}}

\newcommand{\calX}{{\mathcal X}}

\newcommand{\calA}{{\mathcal A}}

\newcommand{\C}{{\mathbb C}}

\newcommand{\R}{{\mathbb R}}

\newcommand{\indexmap}{{\mathrm i}}

\DeclareMathOperator{\tenvec}{vec}

\newcommand{\mb}[1]{\left[\begin{array}{#1}}
\newcommand{\me}{\end{array}\right]}
\newcommand{\smb}{\left[\begin{smallmatrix}}
\newcommand{\sme}{\end{smallmatrix}\right]}

\numberwithin{equation}{section}

\author{Minhong Chen\footnote{Department of Mathematics, Zhejiang Sci-Tech University, Hangzhou, 310029, Zhejiang, P.R.China, {\tt mhchen@zstu.edu.cn}. The work of this author was supported by the National Natural Science Foundation of China (Grant No. 11801513).} \and Daniel Kressner\footnote{Institute of Mathematics, EPF Lausanne, 1015 Lausanne, Switzerland, {\tt daniel.kressner@epfl.ch}.}}

\title{Recursive blocked algorithms for linear systems \\ with Kronecker product structure}

\begin{document}

\maketitle

\begin{abstract}
Recursive blocked algorithms have proven to be highly efficient at the numerical solution of the Sylvester matrix equation and its generalizations. In this work, we show that these algorithms extend in a seamless fashion to higher-dimensional variants of generalized Sylvester matrix equations, as they arise from the discretization of PDEs with separable coefficients or the approximation of certain models in macroeconomics. By combining recursions with a mechanism for merging dimensions, an efficient algorithm is derived that outperforms existing approaches based on Sylvester solvers. 
\end{abstract}

\input{intro}
\input{highdim}
\input{highdim2}

\section{Conclusions, extensions and future work}

We have extended the concept of blocked recursive algorithms to higher-order tensor equations. Both, the complexity estimates and the numerical results, clearly show the importance of combining recursion with merging dimensions in order to arrive at efficient algorithms. For third-order tensor equations, these algorithms seem to constitute the methods of choice. For fourth-order tensor equations with coefficients of nearly equal sizes, reshaping the tensor equation into a Sylvester equation and applying an existing solver is a viable alternative, provided that sufficient memory is available.

The blocked recursive algorithms developed in this work certainly admit extensions to general linear tensor equations taking the form
\[
\sum_{k = 1}^K \bfX \times_1 A_1^{(k)} \times_2 A_2^{(k)}\times_3 \cdots \times_d A_d^{(k)} = \bfB,
\]
assuming that all coefficients $A_\mu^{(k)} \in \R_{n_\mu\times n_\mu}$ are (quasi-)triangular. To transform general coefficients $A_\mu^{(k)}$ into this form requires the existence of invertible matrices $Q_\mu,Z_\mu$ such that $Q_\mu^T A_\mu^{(k)}Z_\mu$ is (quasi-)triangular for every $k = 1,\ldots,K$. For $K \ge 3$, this simultaneous triangularization is only possible under strong additional assumptions on the coefficients. A sufficient condition is that each matrix family $\{ A_\mu^{(1)},\ldots,A_\mu^{(K)} \}$
contains at most two different matrices for $\mu = 1,\ldots,d$. The two classes~\eqref{eq:bigA} and~\eqref{eq:bigAofGsylvTen} appear to constitute the practically most important examples satisfying this condition. 

This work also raises an interesting open question: Is it possible to combine block recursion with low-rank compression, for example in the tensor train format~\cite{Oseledets2011a}, such that the complexity does not grow exponentially with $d$,  assuming that the involved ranks stay constant? It would also be interesting to explore which other numerical linear algebra problems allow for the combination of Kronecker product structure with block recursion. The computation of certain matrix functions, such as the matrix square root~\cite{Deadman2013}, appears to be a likely candidate.

\begin{paragraph}{Acknowledgements.}
Daniel Kressner sincerely thanks Michael Steinlechner and Christine Tobler for insightful discussions on the algorithms presented in this work and their implementation.
\end{paragraph}

\bibliographystyle{plain}
\bibliography{anchp}

\end{document}

%% file: intro.tex

\section{Introduction} \label{sec:intro}

In computations with matrices, recursive blocked algorithms offer an elegant way to arrive at implementations
that benefit from increased data locality and efficiently utilize highly tuned kernels. See~\cite{Elmroth2004}
for a survey and~\cite{Peise2017} for a more recent testimony of this principle. These algorithms
have proven particularly effective for solving Sylvester equations, that is, matrix equations of the form
\begin{equation} \label{eq:trisylv}
 A_1 X + X A_2^T = B,
\end{equation}
where $A_1 \in \R^{n_1\times n_1}, A_2 \in \R^{n_2\times n_2}$, $B \in \R^{n_1\times n_2}$ are given and $X \in \R^{n_1\times n_2}$ is unknown.
In the Bartels-Stewart algorithm~\cite{Bartels1972}, the matrices $A_1$
and $A_2$ are first reduced to block upper form by real Schur decompositions. The reduced
problem is then solved by a variant of backward substitution. Both stages of the algorithms require $O(n^3)$ operations, with $n = \max\{n_1,n_2\}$.
Entirely consisting of level 2 BLAS operations, the backward substitution step performs quite poorly.
To avoid this, Jonsson and K\r{a}gstr\"om~\cite{Jonsson2002,Jonsson2002a} have proposed recursive algorithms for triangular Sylvester and related matrix
equations. The recursive algorithm for solving~\eqref{eq:trisylv} with upper quasi-triangular $A_1,A_2$ starts with partitioning the matrix of larger size.
Assuming $n_1 \ge n_2$, let $A_1 = \begin{pmatrix} A_{1,11} & A_{1,12} \\ 0 & A_{1,22} \end{pmatrix}$
with $A_{1,11} \in \R^{k\times k}$ such that $k\approx n/2$ and partition $X = \begin{pmatrix}  X_1 \\ X_2 \end{pmatrix}$, $B = \begin{pmatrix}  B_1 \\ B_2 \end{pmatrix}$ correspondingly. Then~\eqref{eq:trisylv} becomes equivalent to
\begin{subequations}
 \begin{align} \label{eq:sylv1}
  A_{1,11} X_1 + X_1 A_2^T &= B_1 - A_{1,12} X_2,  \\
 \label{eq:sylv2}
  A_{1,22} X_2 + X_2 A_2^T &= B_2.
 \end{align}
\end{subequations}
First the Sylvester equation~\eqref{eq:sylv2} is solved recursively, then the right-hand side~\eqref{eq:sylv1} is updated, and finally~\eqref{eq:sylv1}
is solved recursively. Apart from the solution of small-sized Sylvester equations at the lowest recursion level, the entire algorithm consists of matrix-matrix multiplications $A_{1,12} X_2$ and thus attains high performance by leveraging level 3 BLAS. As emphasized in~\cite{Elmroth2004,Peise2017}, recursive algorithm are less sensitive to parameter tuning compared to blocked algorithms.

The described algorithm extends to generalized and coupled Sylvester equations, such as $A_1 X M_1 + M_2 X A_2^T = B$; see~\cite{Jonsson2002a,Quintana-Orti2003}. Interestingly, the numerically stable recursive formulation of Hammarling's method~\cite{Hammarling1982} for solving stable Lyapunov equations remains an open problem~\cite{Kressner2008}.

In this paper, we propose several new extensions that address high-dimensional variants of Sylvester equations.
More specifically, we aim at computing a tensor $\bfX \in \R^{n_1\times n_2 \times \cdots \times n_d}$ satisfying the linear equation
\begin{equation} \label{eq:tensorequation}
      \calA \bfX = \bfB,
\end{equation}
where $\calA: \R^{n_1 \times\cdots \times n_d} \to \R^{n_1 \times\cdots \times n_d}$ is a linear operator and $\bfB \in \R^{n_1\times n_2 \cdots \times n_d}$. For $d = 2$, this formulation includes the Sylvester equation~\eqref{eq:trisylv} and its generalizations mentioned above as special cases. For example, for~\eqref{eq:trisylv} the matrix representation of $\calA$ is given by $A = A_2 \otimes I_{n_1} + I_{n_2} \otimes A_1$.

The operator $\calA$ needs to be of a very particular form such that~\eqref{eq:tensorequation} is amenable to the techniques discussed in this work. Motivated by their relevance in applications, we focus on two classes of operators.

\begin{paragraph}{Linear systems with Laplace-like structure.} 
In Section~\ref{sec:laplacelike}, we consider discrete Laplace-like operators $\calA$ having the matrix representation
  \begin{equation} \label{eq:bigA}
         A = A_d \otimes I_{n_{d-1}} \otimes \cdots \otimes I_{n_1} + I_{n_{d}} \otimes A_{d-1} \otimes I_{n_{d-2}} \otimes \cdots \otimes I_{n_1} 
      +\cdots+ I_{n_d} \otimes \cdots \otimes  I_{n_2}  \otimes A_1,
    \end{equation}
with $A_\mu \in \R^{n_\mu \times n_\mu}$, $\mu=1,\ldots,d$. Using the vectorization of tensors,~\eqref{eq:tensorequation} can equivalently be written as $A \tenvec(\bfX) = \tenvec(\bfB)$. Discrete Laplace-like operators arise from the structured discretization of $d$-dimensional PDEs with separable coefficients on tensorized domains. For more general PDEs, matrices of the form~\eqref{eq:bigA} can sometimes be used to construct effective preconditioners; see~\cite{Sangalli2016,Simoncini2016} for examples. Other applications of~\eqref{eq:bigA} arise from Markov chain models~\cite{Dayar2018,Stewart1994b} used, e.g., for simulating interconnected systems. 
\end{paragraph}

\begin{paragraph}{Generalized Sylvester equations with Kronecker structure.} 
Section~\ref{sec:secondclass} is concerned with the second class of operators $\calA$ considered in this work, which have a matrix representation of the form
 \begin{equation} \label{eq:bigAofGsylvTen}
         A = I_{n_d} \otimes I_{n_{d-1}} \otimes \cdots \otimes  I_{n_2}  \otimes A_1 + A_d \otimes A_{d-1} \otimes \cdots \otimes A_2 \otimes C ,
     \end{equation}
with $A_{\mu} \in \R^{n_\mu \times n_\mu}$ for $\mu=1, \ldots, d$ and $C \in \R^{n_1 \times n_1}$. For $d=2$, the linear system~\eqref{eq:tensorequation} now becomes equivalent to the generalized Sylvester equation $A_1 X + C X A_2^T$. For $d>2$, we can view~\eqref{eq:tensorequation} equivalently as a generalized Sylvester equations with coefficients that feature Kronecker structure. If $A_1 = -\lambda I$ for some $\lambda \in \R$ then 
\begin{equation} \label{eq:shiftedkronecker} 
          A = A_d \otimes \cdots \otimes 
          A_2
          \otimes C - \lambda I_{n_1\cdots n_d}.
\end{equation}
Linear systems featuring such shifted Kronecker products have been discussed in~\cite{Martin2006}. The more general case~\eqref{eq:bigAofGsylvTen} arises
from approximations of discrete time DSGE models~\cite{Binning2013}, which play a central role in macroeconomics. 

Recent work on the solution of linear tensor equations~\eqref{eq:tensorequation} has focused on the development of highly efficient approximate and iterative solvers that assume and exploit low-rank tensor structure in the right-hand side and the solution; see~\cite{Grasedyck2013,Hackbusch2012} for overviews. In some cases, these developments can be combined with the methods developed in this work, which do not assume any such structure. For example, if the tensor Krylov subspace method~\cite{Kressner2010} is applied to~\eqref{eq:bigA} for large-scale coefficients $A_\mu$ then our method can be used to solve the smaller-sized linear systems occurring in the method. As far as we know, all existing direct non-iterative solvers for linear tensor equations combine the Bartels-Stewart method for (generalized) Sylvester equations with a recursive traversal of the dimension.
Instances of this approach can be found in~\cite{LiTian2010,Touzene2006} for~\eqref{eq:bigA}, in~\cite{Kamenik2005} for~\eqref{eq:bigAofGsylvTen}, and in~\cite{Martin2006} for~\eqref{eq:shiftedkronecker}. For $d\ge 3$, we are not aware of any work on (recursive) blocked methods that would allow for the effective use of level 3 BLAS.
\end{paragraph}

%
%

%% file: highdim.tex
\section{A recursive blocked algorithm for Laplace-like equations} \label{sec:laplacelike}

Let us first recall two basic operations for tensors from~\cite{Kolda2009}.
The $\mu$th matricization of a tensor $\bfX \in \R^{n_1\times \cdots \times n_d}$ is the matrix $X_{(\mu)} \in \R^{n_\mu\times (n_1\cdots n_{\mu-1}n_{\mu+1}\cdots n_d)}$
obtained by mapping the $\mu$th index to the rows and all other indices to the columns:
\[
   X_{(\mu)}(i_\mu,j ) = \calX(i_1,\ldots,i_d),
\]
with the column index $j$ defined via the index map 
\begin{equation} \label{eq:indexmap}
 j = \indexmap(i_1,\ldots,i_{\mu-1},i_{\mu+1},\ldots,i_d) := 1 + \sum_{\nu = 1\atop \nu\not=\mu}^d (i_\nu -1) \prod_{\eta = 1\atop \eta \not= \mu}^{\nu-1} n_\eta.
\end{equation}
The $\mu$-mode matrix multiplication of $\bfX$ with a matrix $A \in\R^{n_1\times m}$ is the tensor
$\bfY = \bfX \times_\mu A$ satisfying $Y_{(\mu)} = A X_{(\mu)}$. This allows us to rewrite~\eqref{eq:tensorequation}--\eqref{eq:bigA} as
\begin{equation} \label{eq:tensorequation2}
 \bfX \times_1 A_1 + \bfX \times_2 A_2 + \cdots + \bfX \times_d A_d = \bfB.
\end{equation}
It is well known that this equation has a unique solution if and only if $\lambda_1 + \cdots + \lambda_d \not=0$ for any eigenvalues $\lambda_1$ of $A_1$, $\lambda_2$ of $A_2$, etc. In the following, we will assume that this condition is satisfied.

Algorithm~\ref{alg:tensor} describes our general framework for solving~\eqref{eq:tensorequation2}.
Using real Schur decompositions~\cite[Sec. 7.4]{Golub2013}, the coefficient matrices are first transformed to reduced form. More specifically, for each $\mu = 1,\ldots,d$ an orthogonal matrix $U_\mu$ is computed such that $\tilde A_\mu := U_\mu^T A_\mu U_\mu$ is in upper quasi-triangular form, that is, 
$\tilde A_\mu$ is an upper block triangular matrix with $1\times 1$ blocks containing its real eigenvalues and $2\times 2$ blocks containing its complex eigenvalues in conjugate pairs. The right-hand side and the solution tensor need to be transformed accordingly by $\mu$-mode matrix multiplications.
For the rest of this section, we focus on line~\ref{line:reducedeqn} of Algorithm~\ref{alg:tensor}, that is, the solution of the tensor equation with the reduced coefficients.
\begin{algorithm}[H] \small
 \caption{Solution of general Laplace-like equation~\eqref{eq:tensorequation}}
\label{alg:tensor}
\begin{algorithmic}[1]
\STATE Compute real Schur decompositions $A_1 = U_1 \tilde A_1 U_1^T, \ldots, A_d = U_d \tilde A_d U_d^T$.
\STATE Update $\tilde \bfB = \bfB \times_1 U_1^T \times_2 U_2^T \cdots \times_d U^T_d$.
\STATE \label{line:reducedeqn} Compute solution $\tilde \bfX$ of tensor equation~\eqref{eq:tensorequation2} with quasi-triangular coefficients $\tilde A_1,\ldots,\tilde A_d$ and right-hand side $\tilde \bfB$.
\STATE Update $\bfX = \tilde \bfX \times_1 U_1 \times_2 U_2 \cdots \times_d U_d$
\end{algorithmic}
\end{algorithm}

\subsection{Recursion} \label{subsec:highdim_recursion}

By Algorithm~\ref{alg:tensor}, we may assume that $A_1 \in \R^{n_1\times n_1},\ldots,A_d \in \R^{n_d\times n_d}$ are already in upper quasi-triangular form. Choose $\mu$ such that $n_\mu = \max_{\nu}n_\nu$ and $k$ such that $k\approx n_\mu / 2$ and $A_\mu(k+1,k) = 0$.
Partitioning $A_\mu = \begin{pmatrix} A_{\mu,11} & A_{\mu,12} \\ 0 & A_{\mu,22} \end{pmatrix}$ with $A_{\mu,11} \in  \R^{k\times k}$, equation~\eqref{eq:tensorequation2} becomes equivalent to
\begin{subequations}
 \begin{align}
    \bfX_1 \times_\mu A_{\mu,11} + \sum_{\nu = 1 \atop \nu\not=\mu}^d \bfX_1 \times_\nu A_\nu  &= \bfB_1 - \bfX_2 \times_\mu A_{\mu,12}, \label{eq:tensorequation11} \\
  \bfX_2 \times_\mu A_{\mu,22} + \sum_{\nu = 1 \atop \nu\not=\mu}^d \bfX_2 \times_\nu A_\nu  &= \bfB_2, \label{eq:tensorequation22}
 \end{align}
\end{subequations}
where \begin{subequations} \label{eq:splitX}
 \begin{align}
\bfX_1 & = \bfX(1:n_1,\ldots,1:n_{\mu-1},1:k,1:n_{\mu+1},\ldots,1:n_{d}), \\
\bfX_2 & = \bfX(1:n_1,\ldots,1:n_{\mu-1},k+1:n_\mu,1:n_{\mu+1},\ldots,1:n_{d}),
\end{align}
\end{subequations}
and $\bfB_1,\bfB_2$ are defined analogously. Noting that~\eqref{eq:tensorequation22} and~\eqref{eq:tensorequation11} are again equations with Laplace-like operators, they can be solved recursively. The recursion is stopped once the maximal size is below a user-specified block size $n_{\min} \ge 2$.
These considerations lead to Algorithm~\ref{alg:laplacerecursive}.
\begin{algorithm}[H] \small
 \caption{Recursive sol. of quasi-triang. tensor eqn: $\bfX=\mathtt{reclap}(A_1,\ldots,A_d,\bfB)$ }
\label{alg:laplacerecursive}
\begin{algorithmic}[1]
\STATE Determine $n_\mu = \max_{\nu}n_\nu$.
\STATE \label{line:full} {\bf if} $n_\mu \le n_{\min}$ {\bf then} solve full system and return {\bf end if}
\STATE Set $k = \lfloor n_\mu/2 \rfloor$. {\bf if} $A_\mu(k+1,k)\not=0$ {\bf then} $k = k + 1$ {\bf end if}
\STATE Set $\bfB_1 = \bfB(:,\ldots,:,1:k,:,\ldots,:)$, $\bfB_2 = \bfB(:,\ldots,:,k+1:n_\mu,:,\ldots,:)$.
\STATE Call $\bfX_2=\mathtt{reclap}(A_1,\ldots,A_{\mu-1},A_{\mu}(k+1:n_\mu,k+1:n_\mu),A_{\mu+1},\ldots,A_d,\bfB_2)$.
\STATE \label{line:matmatrec} Update $\bfB_1 \leftarrow \bfB_1 - \bfX_2 \times_\mu A_{\mu}(1:k,k+1:n_\mu)$.
\STATE Call $\bfX_1=\mathtt{reclap}(A_1,\ldots,A_{\mu-1},A_{\mu}(1:k,1:k),A_{\mu+1},\ldots,A_d,\bfB_2)$.
\STATE Concatenate $\bfX_1, \bfX_2$ along $\mu$th mode into $\bfX$.
\end{algorithmic}
\end{algorithm}

Let $\mathsf{comp}(n)$ denote the complexity of Algorithm~\ref{alg:laplacerecursive} for even $n = n_1 = \cdots = n_d$. 
On the top level of recursion Algorithm~\ref{alg:laplacerecursive} is applied to one $n\times \cdots \times n$ tensor, on the second level to two 
$n/2 \times n \times \cdots \times n$ tensors, on the third level to four 
$n/2 \times n/2 \times n \times \cdots \times n$ tensors, and so on. Under the slightly simplified assumption that the multiplication of an $n/2 \times n/2$ quasi-triangular matrix with a vector requires $n^2 / 4$ floating point operations (flops), each level of the first $d$ recursions requires a total of $n^{d+1}/4$ flops to execute the matrix-matrix multiplications in line~\ref{line:matmatrec} of Algorithm~\ref{alg:laplacerecursive}. After $d$ recursions of Algorithm~\ref{alg:laplacerecursive}, $n$ has been reduced to $n/2$ in each mode and, therefore,
$
 \mathsf{comp}(n) = d n^{d+1} / 4 + 2^d  \mathsf{comp}(n/2).
$ Assuming that $n/n_{\min}$ is a power of two, we obtain
\begin{equation} \label{eq:costreductiontonmin}
 \mathsf{comp}(n) = O\big( n^{d+1} \big)  + (2^d)^{\log_2 n/n_{\min}} 
 \mathsf{comp}(n_{\min}) = O\big( n^{d+1} \big)  + \frac{n^d}{n^d_{\min}}   
 \mathsf{comp}(n_{\min}).
\end{equation}
Once the maximal size of the tensor is $n_{\min}$ or below, line~\ref{line:full} of Algorithm~\ref{alg:laplacerecursive} assembles the 
matrix $A$ defined in~\eqref{eq:bigA} and solves the block triangular linear system $A \tenvec(\bfX) = \tenvec(\bfB)$ by backward substitution.
This requires $O\big( (n_{\min})^{2d}\big)$ flops and therefore  
\begin{equation} \label{eq:recursivecomplexity}
  \mathsf{comp}(n) = O\big( n^{d+1} + n^d_{\min} n^{d} \big).
\end{equation}
This compares favorably with the $O(n^{2d})$ operations needed by backward substitution applied to the assembled full triangular linear system.
The complexity estimate~\eqref{eq:recursivecomplexity} also reflects the critical role played by the solution of the small systems in line~\ref{line:full}.
On the one hand, the operation count suggests to choose $n_{\min}$ as small as possible, say, $n_{\min} = 2$. On the other hand, it has been observed for $d = 2$ in~\cite{Jonsson2002} that a small value of $n_{\min}$ creates significant overhead and requires very well tuned kernels. In the following section, we describe a technique that alleviates this difficulty.

\subsection{Merging dimensions: triangular case}
\label{subsec:merge_highdim}

To avoid the critical dependence on $n_{\min}$ observed in~\eqref{eq:recursivecomplexity} we replace line~\ref{line:full} of Algorithm~\ref{alg:laplacerecursive} by the following procedure. Once $n_1 n_2 \le n_{\min}^2$, the matrix
\begin{equation}
    \label{eq:merged_equation}
 A^\prime_{1} = I_{n_2} \otimes A_1 + A_2 \otimes I_{n_1} 
\end{equation}
is formed explicitly. For the moment, let us suppose that $A_1$ and $A_2$ are upper triangular. This can be achieved by computing complex instead of real Schur decompositions in Algorithm~\ref{alg:tensor}, leading to a triangular tensor equation with complex coefficients. Because of roundoff error, the computed solution to the original equation will now have a (small) imaginary part. This can be safely set to zero~\cite{Martin2006b}.

The matrix $A^\prime_{1}$ inherits the triangular structure from $A_1,A_2$ and the $d$-dimensional tensor equation~\eqref{eq:tensorequation2} is equivalent to the $(d-1)$-dimensional equation
\begin{equation} \label{eq:mergedm1}
 \bfX^\prime \times_1 A^\prime_{1} + \bfX^\prime \times_3 A_{3} + \cdots + \bfX^\prime \times_d A_{d} = \bfB^\prime,
\end{equation}
with reshaped $\bfX^\prime, \bfB^\prime \in \C^{n_1n_2 \times n_3 \cdots \times n_d}$. This equation is solved recursively. A major advantage, this approach allows us to reduce $d$. For $d=3$, the system~\eqref{eq:mergedm1} becomes the triangular Sylvester equation
\[
 A^\prime_1 \bfX^\prime + \bfX^\prime A_3^T = \bfB^\prime,
\]
to which the efficient solvers described in Section~\ref{sec:intro} can be applied. Note that $A_3^T$ now refers to the \emph{complex} transpose of $A_3 \in \C^{n_3\times n_3}$. Algorithm~\ref{alg:laplacemerge} summarizes the proposed procedure.

{\small 
\begin{algorithm}[H] \small
\caption{Recursive sol. of triangular tensor equation: $\bfX=\mathtt{mrglap}(A_1,\ldots,A_d,\bfB)$\label{alg:laplacemerge}}
\begin{algorithmic}[1]
\IF{$n_1 n_2 \le n_{\min}^2$} \label{line:ifcondition}
 \STATE Compute $A^\prime_{1} = I_{n_2} \otimes A_1 + A_2 \otimes I_{n_1}$. \label{line:merge}
 \STATE Reshape $\bfB^\prime(1:n_1 n_2,:,\ldots,:) = \bfB(1:n_1,1:n_2,:,\ldots,:)$.
 \IF{$d = 3$}
 \STATE Solve Sylvester equation $A^\prime_1 \bfX^\prime + \bfX^\prime A_3^T = \bfB^\prime$. \label{line:sylvester}
 \ELSE
 \STATE Call $\bfX^\prime=\mathtt{mrglap}(A_1^\prime,A_3,\ldots,A_d,\bfB^\prime)$
 \ENDIF
 \STATE Reshape $\bfX(1:n_1,1:n_2,:,\ldots,:) = \bfX^\prime(1:n_1 n_2,:,\ldots,:)$.
\ELSE 
\STATE Determine $n_\mu = \max_{\nu}n_\nu$ and set $k = \lfloor n_\mu/2 \rfloor$. 
\STATE Set $\bfB_1 = \bfB(:,\ldots,:,1:k,:,\ldots,:)$, $\bfB_2 = \bfB(:,\ldots,:,k+1:n_\mu,:,\ldots,:)$.
\STATE Call $\bfX_2=\mathtt{mrglap}(A_1,\ldots,A_{\mu-1},A_{\mu}(k+1:n_\mu,k+1:n_\mu),A_{\mu+1},\ldots,A_d,\bfB_2)$.
\STATE \label{line:matmat} Update $\bfB_1 \leftarrow \bfB_1 - \bfX_2 \times_\mu A_{\mu}(1:k,k+1:n_\mu)$.
\STATE Call $\bfX_1=\mathtt{mrglap}(A_1,\ldots,A_{\mu-1},A_{\mu}(1:k,1:k),A_{\mu+1},\ldots,A_d,\bfB_2)$.
\STATE Concatenate $\bfX_1, \bfX_2$ along $\mu$th mode into $\bfX$.
\ENDIF
\end{algorithmic}
\end{algorithm}
}

To analyze the complexity of Algorithm~\ref{alg:laplacemerge} for $n_1 = \cdots = n_d = n > 2 n_{\min}$, we observe that all sizes are first reduced to $2 n_{\min}$ or below before the condition in line~\ref{line:ifcondition} is met. Hence, up to constant factors the recursive estimate~\eqref{eq:costreductiontonmin} holds and it remains to discuss the complexity for $n_1 = \cdots = n_d = n_{\min}$, which will be denoted by $\overline{\mathsf{comp}}_d(n_{\min})$. The merge in line~\ref{line:merge} reduces the order to $d-1$ but increases the first mode size to $n_{\min}^2$. Approximately $\log_2 ( n_{\min}^2 / n_{\min} ) = \log_2 n_{\min}$ recursions are needed to reduce it back to $n_{\min}$. Similarly as in Section~\ref{subsec:highdim_recursion} we calculate
\[
 \overline{\mathsf{comp}}_d(n_{\min}) = O\big( n_{\min}^{d+2} \big) + n_{\min} \overline{\mathsf{comp}}_{d-1}(n_{\min}) = 
 O\big( n_{\min}^{d+2} \big) + n^{d-3}_{\min}\, \overline{\mathsf{comp}}_{3}(n_{\min}).
\]
For $d =3$, the solution of the triangular Sylvester equation in line~\ref{line:sylvester} requires $O\big( n_{\min}^{5} \big)$ flops. In turn,
$\overline{\mathsf{comp}}_d(n_{\min}) = O\big( n_{\min}^{d+2} \big)$.
Inserted into~\eqref{eq:costreductiontonmin}, we arrive at 
\[
O\big(n^{d+1} + n_{\min}^{2} n^{d}\big)
\]
flops for Algorithm~\ref{alg:laplacemerge}. For $d \ge 3$, this compares favorably with the complexity estimate~\eqref{eq:recursivecomplexity} for Algorithm~\ref{alg:laplacerecursive}; the dependence on $n_{\min}$ has been reduced significantly. Equally importantly, Algorithm~\ref{alg:laplacemerge} allows us to leverage efficient solvers for triangular Sylvester equations, such as the ones described in~\cite{Jonsson2002}.

\subsection{Merging dimensions: quasi-triangular case} \label{sec:quasitriangular}

The use of complex arithmetic, which increases the cost (by a constant factor) in terms of operations and memory, can be avoided when using the real Schur form and working with quasi-triangular coefficients. However, a few modifications are needed because the matrix $A^\prime_{1}$ formed in~\eqref{eq:merged_equation} does \emph{not} inherit the quasi-triangular structure from $A_1$ and $A_2$. To illustrate what happens, let us consider the following example for $n_1 = 3, n_2 = 4$:
\begin{equation} \label{eq:a1}
A_{1}={\scriptsize\left[ \arraycolsep=1.4pt\def\arraystretch{1.1}
\begin{array}{cccc}   \times &  \times & \times &  \times \\   0 & \times & \times &  \times \\   0 & \times & \times & \times \\    0 & 0 & 0 & \times \end{array} \right]}, \quad
 A_2={\scriptsize\left[ \arraycolsep=1.4pt\def\arraystretch{1.1} \begin{array}{ccc} \times &  \times & \times \\  0 & \times & \times \\  0 & \times & \times \end{array} \right]},\quad
 A_1^\prime = {\scriptsize 
\left[ \arraycolsep=1.4pt\def\arraystretch{1.1}
\begin{array}{cccc|cccc|cccc}
\times &  \times & \times &  \times & \times & 0 & 0 & 0 & \times & 0 & 0 & 0  \\ 
0 &  \times & \times &  \times & 0 & \times & 0 & 0 & 0 & \times & 0 & 0  \\ 
0 &  \times & \times &  \times & 0 & 0 & \times & 0 & 0 & 0 & \times & 0  \\ 
0 & 0 & 0 &  \times & 0 & 0 & 0 & \times & 0 & 0 & 0 & \times \\ \hline
0 & 0 & 0 & 0 & \times & \times& \times& \times & \times & 0 & 0 & 0 \\
0 & 0 & 0 & 0 & 0 & \times& \times& \times & 0 & \times& 0 & 0 \\
0 & 0 & 0 & 0 & 0 & \times & \times & \times& 0 & 0 &\times & 0 \\
0 & 0 & 0 & 0 & 0 & 0 & 0 & \times& 0& 0 & 0 & \times \\ \hline 
0 & 0 & 0 & 0 & \times & 0 & 0 & 0 & \times & \times& \times& \times  \\
0 & 0 & 0 & 0 & 0 & \times & 0 & 0 & 0 & \times& \times& \times  \\
0 & 0 & 0 & 0 & 0 & 0 & \times & 0 & 0 & \times & \times & \times  \\
0 & 0 & 0 & 0 & 0 & 0 & 0 & \times & 0 & 0 & 0 & \times 
\end{array}
\right]}.
\end{equation}
The diagonal matrix at the (3,2) block disturbs the quasi-triangular structure of $A_1^\prime$. More generally, assuming $n_1 = n_2 = n_{\min}$ the matrix $A_1^\prime$ is an $n_{\min}^2 \times n_{\min}^2$ block upper triangular matrix with diagonal blocks of size at most $n_{\min}$. This matrix can be returned to quasi-triangular form by computing a real Schur decomposition of $A_1^\prime$. The impact of this operation on the overall cost of Algorithm~\ref{alg:laplacemerge} can be made negligible by exploiting the structure of $A_1^\prime$:

\begin{itemize}
 \item When the structure of $A_1^\prime$ is completely ignored, its real Schur decomposition takes $O(n_{\min}^6)$ flops and, in turn, the complexity of Algorithm~\ref{alg:laplacemerge} increases to $O\big(n^{d+1} + n_{\min}^{3} n^{d}\big)$. 
 \item When the block triangular structure of $A_1^\prime$ is taken into account, the cost of computing its real Schur decomposition reduces to $O(n_{\min}^5)$ flops. When used within Algorithm~\ref{alg:laplacemerge}, the additional flops spent on performing these decompositions and applying the resulting orthogonal transformations amounts to 
 $O\big( n_{\min}^{2} n^{d} \big)$ in total. In turn, this operation does not increase the complexity of Algorithm~\ref{alg:laplacemerge} but its dependence on $n_{\min}^{2}$ is not negligible either.
 \item The diagonal structure of the off-diagonal blocks of $A_1^\prime$ can be exploited to reduce the cost further, using a permutation trick similar to the one discussed in~\cite{Martin2006}. To illustrate this, consider the $12\times 12$ matrix $A_1^\prime$ from~\eqref{eq:a1}. By applying a perfect shuffle permutation~\cite{VanLoan2000} to the last $8$ rows and columns,
 we obtain the permuted matrix
 \[
  P^T A_1^\prime P = {\scriptsize 
\left[ \arraycolsep=1.4pt\def\arraystretch{1.1}
\begin{array}{cccc|cc|cccc|cc}
\times &  \times & \times &  \times & \times &  \times & 0 & 0 & 0 & 0 & 0 & 0  \\ 
0 & \times & \times &  \times & 0 & 0 & \times & \times& 0 & 0 &0 & 0\\ 
0 &  \times & \times &  \times & 0 & 0 & 0 & 0 & \times & \times & 0 & 0  \\ 
0 & 0 & 0 &  \times & 0 & 0 & 0 & 0 & 0 & 0 & \times & \times \\ \hline
0 & 0 & 0 & 0 & \times & \times& \times& 0 & \times & 0 & \times & 0 \\
0 & 0 & 0 & 0 & \times & \times& 0& \times & 0 & \times& 0 & \times \\ \hline
0 & 0 & 0 & 0 & 0 & 0 & \times & \times& \times & 0 &\times & 0 \\
0 & 0 & 0 & 0 & 0 & 0 & \times & \times& 0& \times & 0 & \times \\ 
0 & 0 & 0 & 0 & 0 & 0 & \times & 0 & \times & \times& \times& 0 \\
0 & 0 & 0 & 0 & 0 & 0 & 0 & \times & \times & \times& 0 & \times  \\ \hline 
0 & 0 & 0 & 0 & 0 & 0 & 0 & 0 & 0 & 0 & \times & \times  \\
0 & 0 & 0 & 0 & 0 & 0 & 0 & 0 & 0 & 0 & \times & \times 
\end{array}
\right].}
 \]
 In the general case, applying such a permutation to each $n_{\min} \times n_{\min}$
 diagonal block transforms $A_1^\prime$ into a block upper triangular matrix with diagonal blocks of size at most $4$. This reduces the cost of  computing its real Schur decomposition to $O(n_{\min}^4)$ flops and the overall impact of this operation on the cost of Algorithm~\ref{alg:laplacemerge} becomes negligible.
\end{itemize}

\subsection{Numerical experiments} \label{sec:numexplaplace}

All algorithms proposed in this work have been implemented in {\sc Matlab} R2019a and executed on a Lenovo ThinkPad T460, which comes with an Intel Core i5-6300U processor and 8 Gbytes of DDR3L-RAM. The implementation of the algorithms together with scripts for reproducing each of the experiments reported in this work are available from \url{https://anchp.epfl.ch/misc/}.

Care has been taken to avoid unnecessary overhead in our {\sc Matlab} implementation.
For example, the {\tt tensor} object from the Tensor Toolbox \cite{Bader2012} is very convenient for realizing tensor operations but our preliminary experiments indicated that its use in Algorithms~\ref{alg:laplacerecursive} and~\ref{alg:laplacemerge} would lead to significant performance loss, possibly due to excessive memory transfer. Instead, we directly use {\sc Matlab} arrays, combined with the {\tt permute} and {\tt reshape} functions for implementing $\mu$-mode matrix multiplications. For solving triangular Sylvester equations, as needed, e.g., in Algorithm~\ref{alg:laplacemerge}, we utilize the internal {\sc Matlab} function {\tt sylvester\_tri}. This function seems to be based on the algorithms presented in~\cite{Jonsson2002,Jonsson2002a} and avoids performing any additional Schur decomposition. 

The techniques from Section~\ref{sec:quasitriangular}, which allow for the use of real arithmetic in Algorithm~\ref{alg:laplacemerge}, have been implemented and verified. However, we observed that none of the three described variants leads to competitive performance, any benefit from structure exploitation is offset by the overhead it incurs in {\sc Matlab}, due to the relatively small values of $n_{\min}$ needed for reaching good performance. In the following, we therefore consistently use complex Schur decompositions for reducing all coefficients to triangular form. All reported times include the time needed by Algorithm~\ref{alg:tensor} for performing these decompositions and applying the corresponding transformations. The coefficients used in our experiments have been generated with {\tt randn}. 

\begin{figure}
    \centering
    \subfloat[Algorithm~\ref{alg:laplacerecursive}]{\includegraphics[width=0.49\textwidth]{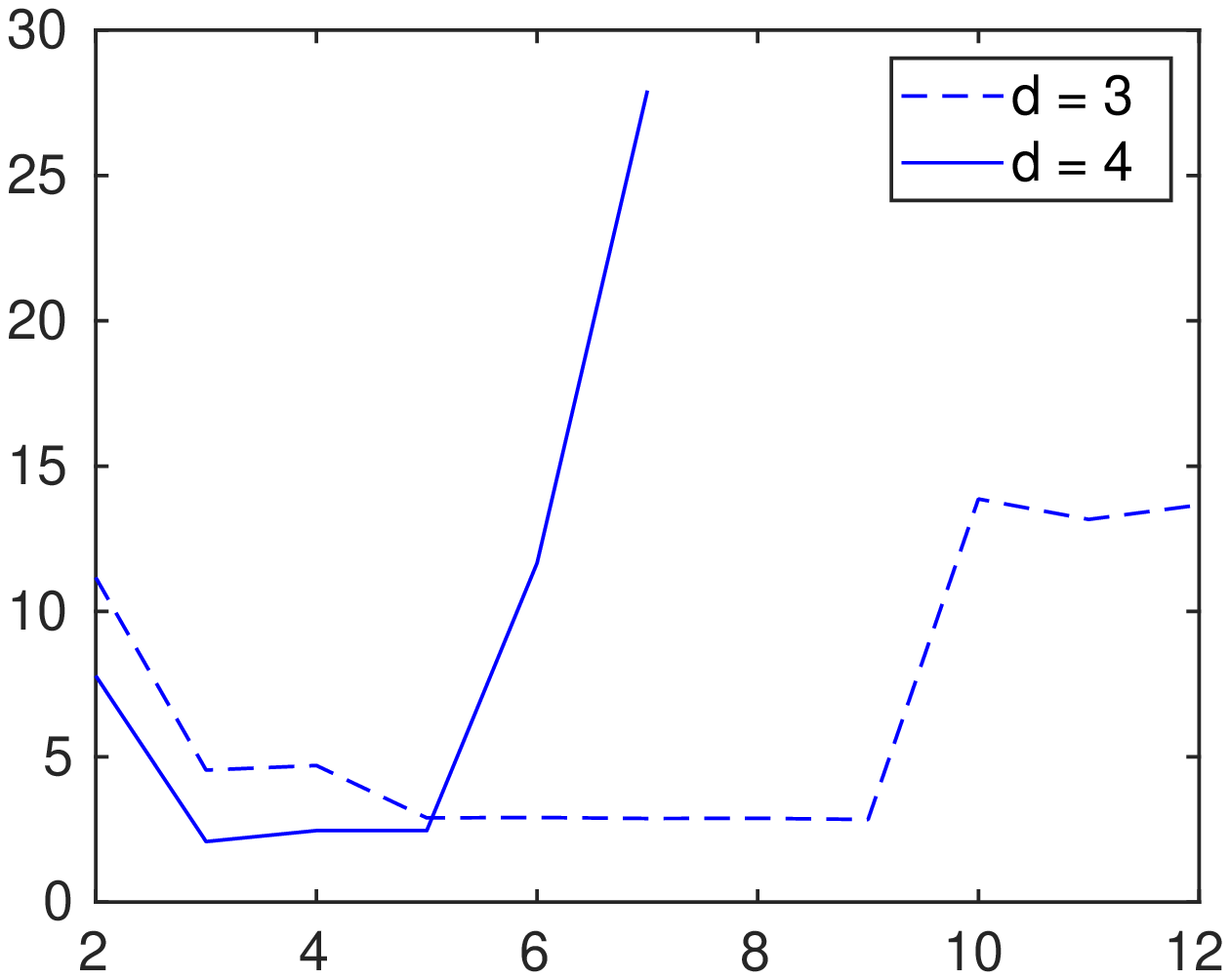}} 
    \subfloat[Algorithm~\ref{alg:laplacemerge}]{\includegraphics[width=0.49\textwidth]{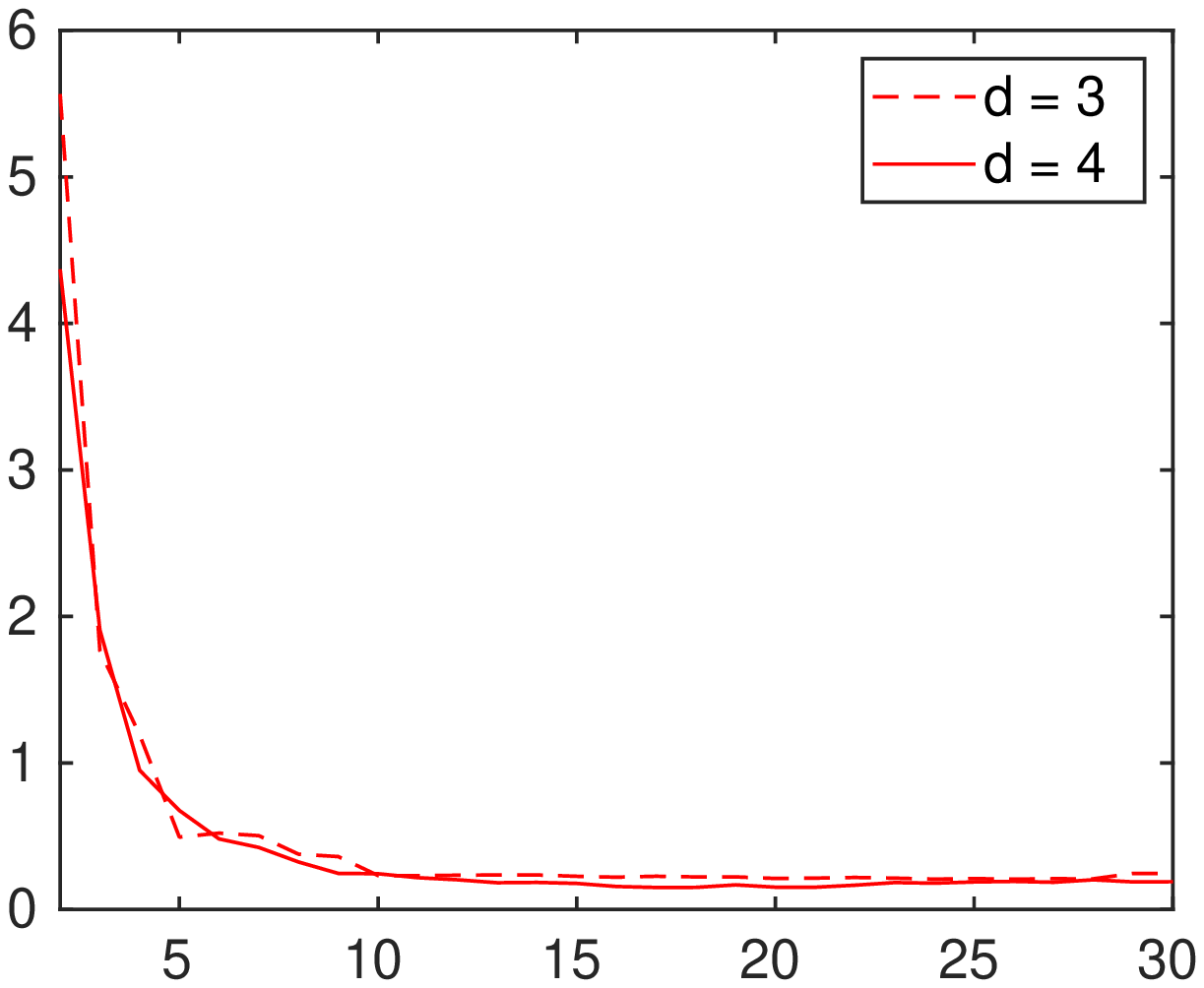}} 
    \caption{Execution times (in seconds) vs. $n_{\min}$ for Algorithms~\ref{alg:laplacerecursive} and~\ref{alg:laplacemerge} applied to random $n\times\cdots\times n$ tensors with $n = 80$ for $d = 3$ and $n = 25$ for $d = 4$. \label{fig:nminlaplace}}
\end{figure}

\begin{paragraph}{Choice of $n_{\min}$.} Figure~\ref{fig:nminlaplace} shows the execution times obtained for fixed $n$ and varying $n_{\min}$. All numbers have been averaged over five consecutive runs. As to be expected from the complexity estimates, the performance of Algorithm~\ref{alg:laplacerecursive} is very sensitive to the choice of $n_{\min}$, especially for $d = 4$. The smallest execution times are attained by $n_{\min} = 7$ for $d = 3$ and $n_{\min} = 3$ for $d = 4$. 
The performance of Algorithm~\ref{alg:laplacemerge} is not very sensitive to the choice of $n_{\min}$, provided that its value is not chosen too small.
The smallest execution times are attained by $n_{\min} = 26$ for $d = 3$, $n_{\min} = 18$ for $d = 4$, and $n_{\min} = 14$ for $d = 5$. These values of $n_{\min}$ are used in the following.
\end{paragraph}

\begin{figure}
\begin{minipage}{.5\linewidth}
\centering
\subfloat[$d = 3$]{\includegraphics[scale=.5]{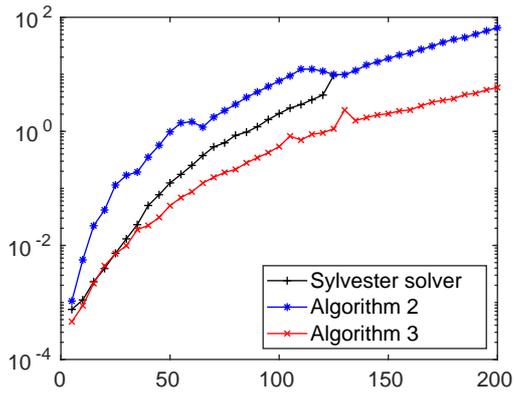}}
\end{minipage}%
\begin{minipage}{.5\linewidth}
\centering
\subfloat[$d = 4$]{\includegraphics[scale=.5]{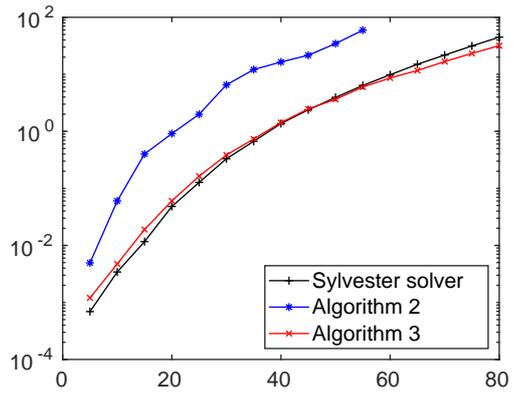}}
\end{minipage} \\
\centering
\subfloat[$d = 5$]{\includegraphics[scale=.5]{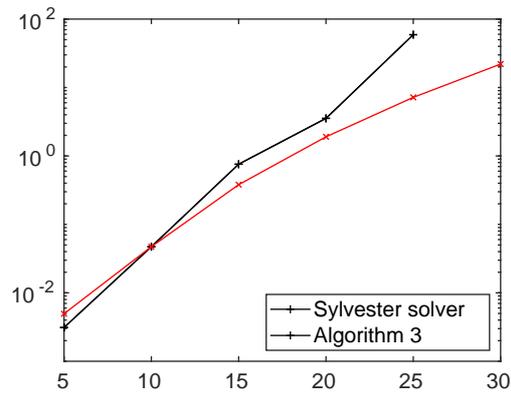}}
    \caption{Execution times (in seconds) vs. $n$ for Algorithms~\ref{alg:laplacerecursive} and~\ref{alg:laplacemerge} compared to Sylvester solver. \label{fig:complaplace}}
\end{figure}

\begin{paragraph}{Comparison.}
We have compared our newly proposed algorithms with the following procedure termed ``Sylvester solver'': After reducing the coefficients $A_1,\ldots, A_d$ of the Laplace-like equation~\eqref{eq:tensorequation2} to triangular form and reshaping $\bfB$ suitably into a matrix $B$, one of the Sylvester equations 
\[
\begin{array}{l}
  (I \otimes A_1  + A_2 \otimes I) X + X A_3^T = B, \quad 
 (I \otimes A_1  + A_2 \otimes I) X + X (I \otimes A_3  + A_4 \otimes I)^T = B \\
 (I \otimes A_1  + A_2 \otimes I) X + X (I \otimes I \otimes A_3  + I\otimes A_4 \otimes I + A_5\otimes I \otimes I)^T = B
 \end{array}
\]
is solved for $d = 3, 4, 5$ by calling {\tt sylvester\_tri}. The results reported in Figure~\ref{fig:complaplace} confirm that Algorithms~\ref{alg:laplacerecursive} and~\ref{alg:laplacemerge} have the same asymptotic cost. However, Algorithm~\ref{alg:laplacemerge} is always faster, by an order of magnitude for sufficiently large $n$.  
For $d = 3$, the Sylvester solver is nearly always slower than Algorithm~\ref{alg:laplacemerge}. For $d = 4$, the picture is less clear; only for  $n\ge50$ becomes  Algorithm~\ref{alg:laplacemerge}, which has complexity $O(n^5)$, consistently faster than the Sylvester solver, which has complexity $O(n^6)$. For $d=5$, the difference in complexity is more pronounced and, in turn, Algorithm~\ref{alg:laplacemerge} is nearly always faster.
\end{paragraph}

For all experiments performed, the norm of the residual was checked and no significant differences in terms of numerical stability were observed between the different algorithms tested.

%% file: highdim2.tex
\section{A recursive blocked algorithm for generalized Sylvester equations with Kronecker structure} \label{sec:secondclass}

In this section, we extend the developments from Section~\ref{sec:laplacelike} to the second class of operators $\calA$ considered in this work, which have the matrix representation~\eqref{eq:bigAofGsylvTen}.
The corresponding linear system reads in tensor notation as
\begin{equation} \label{eq:tensorequation3}
 \bfX \times_1 A_1 + \bfX \times_1 C \times_2 A_2 \times_3 A_3 \cdots \times_d A_d = \bfB.
\end{equation}
Because of its connection to generalized Sylvester equations~\cite{Chu1987} explained in the introduction, this equation has a unique solution if and only if the matrix pencil $A_1+\lambda C$ is regular and none of its eigenvalues is an eigenvalue of $-A_d \otimes A_{d-1} \otimes \cdots \otimes A_2$. In the following, we assume that this condition is satisfied.

Algorithm~\ref{alg:gstensor} is the equivalent of Algorithm~\ref{alg:tensor} for reducing~\eqref{eq:tensorequation3} to quasi-triangular form. The most notable difference is that now a generalized Schur decomposition~\cite[Sec. 7.7.2]{Golub2013} of $A_1+\lambda C$ needs to be computed, using the QZ algorithm.

\begin{algorithm}[H] \small
 \caption{Solution of generalized Sylvester equation~\eqref{eq:tensorequation3}}
\label{alg:gstensor}
\begin{algorithmic}[1]
\STATE Compute real generalized Schur decomposition $A_1 = U_1 \tilde A_1 Z^T$,  $C=U_1 \tilde C Z^T$.
\STATE Compute real Schur decompositions $A_2 = U_2 \tilde A_2 U_2^T,  \cdots, A_d = U_d \tilde A_d U_d^T$.
\STATE Update $\tilde \bfB = \bfB \times_1 U_1^T \times_2 U_2^T \cdots \times_d U^T_d$.
\STATE \label{line:reducedeq3} Compute solution $\tilde \bfX$ of tensor equation~\eqref{eq:tensorequation3} with (quasi-)triangular coefficients $\tilde C,\tilde A_1,\ldots,\tilde A_d$ and right-hand side $\tilde \bfB$.
\STATE Update $\bfX =  \tilde{\bfX} \times_1 Z \times_2 U_2 \cdots \times_d U_d$
\end{algorithmic}
\end{algorithm}

\subsection{Recursion} \label{subsec:highdim2_recursion}

The rest of this section is concerned with line~\ref{line:reducedeq3} of Algorithm~\ref{alg:gstensor}, solving~\eqref{eq:tensorequation3} with  
upper quasi-triangular coefficients $A_1 \in \R^{n_1\times n_1},\ldots,A_d \in \R^{n_d\times n_d}$ and upper triangular $C \in \R^{n_1\times n_1}$.
Again we proceed recursively and choose $\mu$ such that $n_\mu = \max_{\nu}n_\nu$ and $k$ such that $k\approx n_\mu / 2$ and $A_\mu(k+1,k) = 0$. We partition $A_\mu = 
\begin{pmatrix} A_{\mu,11} & A_{\mu,12} \\ 0 & A_{\mu,22} \end{pmatrix}$, $A_{\mu,11} \in \R^{k\times k}$ and split
the tensors $\bfX$ and $\bfB$ along their $\mu$th mode into $\bfX_1,\bfX_2$ and $\bfB_1,\bfB_2$, respectively, in accordance with~\eqref{eq:splitX}.

\textbf{Case 1: $\mu=1$.} We additionally partition $C = 
\begin{pmatrix} C_{11} & C_{12} \\ 0 & C_{22} \end{pmatrix}$ and decouple~\eqref{eq:tensorequation3} along the first mode:
 \begin{align*}  
\bfX_1 \times_1 A_{1,11} + \bfX_1 \times_1 C_{11}\times_2 A_2 \times_3 \cdots \times_d A_d  
&= \hat \bfB_1,  \\
\bfX_2 \times_1 A_{1,22} +\bfX_2 \times_1 C_{22} \times_2 A_2 \times_3 \cdots \times_d A_d & = \bfB_2,
 \end{align*}
with $\hat \bfB_1:= \bfB_1 - \bfX_2 \times_1 A_{1,12} -\bfX_2 \times_1 C_{12}\times_2 A_2 \times_3 \cdots \times_d A_d$.
Both equations take the form of the tensor equation~\eqref{eq:tensorequation3} with (quasi-)triangular coefficients.  We recursively solve for $\bfX_2$ and then solve for $\bfX_1$, after computing $\hat \bfB_1$. 

\textbf{Case 2: $\mu \neq1$.} Decoupling~\eqref{eq:tensorequation3} along the $\mu$th mode gives the two tensor equations
 \begin{align*}  
\bfX_1 \times_1 A_1 + \bfX_1 \times_1 C \times_2  \cdots \times_\mu A_{\mu,11} \times_{\mu+1} \cdots \times_d A_d &= \hat \bfB_1, \\
\bfX_2 \times_1 A_1 +\bfX_2 \times_1  C \times_2  \cdots \times_\mu A_{\mu,22} \times_{\mu+1} \cdots \times_d A_d &=\bfB_2,
\end{align*}
with $\hat \bfB_1:= \bfB_1 - \bfX_2  \times_1 C \times_2  \cdots \times_\mu A_{\mu,12} \times_{\mu+1} \cdots \times_d A_d$.
Again, we first solve for $\bfX_2$ and then for $\bfX_1$.

\begin{algorithm}[H] \small
 \caption{Recursive solution of triangular generalized Sylvester equation with tensor structure: $\bfX=\mathtt{recgsylvten} (A_1,\ldots, A_d, C, \bfB)$}
\label{alg:gsylvtenrec}
\begin{algorithmic}[1]
\STATE  Determine $n_\mu = \max_{\nu}n_\nu$.
\STATE  \label{line:gsylvfull} \textbf{if} $n_\mu \leq n_{\min}$ \textbf{then} solve full system and return \textbf{end if}
\STATE  Set $k = \lfloor n/2 \rfloor$. {\bf if} $A_\mu(k+1,k)\not=0$ {\bf then} $k = k + 1$ {\bf end if}
\STATE \label{line:startcode} Set $\bfB_1 = \bfB(:,\ldots,:,1:k,:,\ldots,:)$, $\bfB_2 = \bfB(:,\ldots,:,k+1:n_\mu,:,\ldots,:)$.
   \IF {  $\mu=1$}        
            \STATE Call $\bfX_2=\mathtt{recgsylvten}(A_1(k+1:n_1,k+1:n_1), A_2, \ldots, A_d, C(k+1:n_1, k+1:n_1), \bfB_2)$.
	    \STATE \label{line:matmat1} Update $ \bfB_1 \leftarrow  \bfB_1 - \bfX_2  \times_1A_1(1:k, k+1:n_1) - \bfX_2  \times_1 C(1:k, k+1:n_1)   \times_2 A_2 \times_3 \ldots \times_d  A_d$.
	    \STATE Call $\bfX_1=\mathtt{recgsylvten}(A_1(1:k, 1:k), A_2  \ldots,A_d, C(1:k, 1:k), \bfB_2)$.
\ELSE
\STATE Call $\bfX_2=\mathtt{recgsylvten}(A_1,\ldots,A_{\mu-1},A_{\mu}(k+1:n_\mu,k+1:n_\mu),A_{\mu+1},\ldots, A_d, C, \bfB_2)$.
\STATE  Update \label{line:matmat2} $ \bfB_1 \leftarrow  \bfB_1 -  \bfX_2  \times_1 C  \times_2 A_2 \times_3 \cdots \times_\mu A_{\mu}(1:k, k+1:n_\mu)\times_{\mu+1}A_{\mu+1} \times_{\mu+2} \ldots \times_d A_d$.
\STATE Call $\bfX_1=\mathtt{recgsylvten}(A_1,\ldots,A_{\mu-1},A_{\mu}(1:k,1:k),A_{\mu+1},\ldots, A_d, C, \bfB_2)$.
\ENDIF
\STATE \label{line:endcode} Concatenate $\bfX_1, \bfX_2$ along $\mu$th mode into $\bfX$.
\end{algorithmic}
\end{algorithm}

Algorithm~\ref{alg:gsylvtenrec} summarizes the described procedure. Compared to Algorithm~\ref{alg:laplacerecursive}, the largest difference is that the right-hand side updates in lines~\ref{line:matmat1} and~\ref{line:matmat2} require up to $d$ matrix multiplications instead of only one. While potentially having an impact on computational time, this has no impact on the asymptotic complexity, which remains $O\big( n^{d+1} + n_{\min}^d n^d \big)$.


\subsection{Merging dimensions: triangular case}
\label{subsec:merge_highdim2}

In analogy to the discussion in Section~\ref{subsec:merge_highdim}, we now discuss the combination of Algorithm~\ref{alg:gsylvtenrec} with a merging procedure that helps to alleviate the critical dependence of its performance on $n_{\min}$. Again, we first suppose that all coefficients triangular. This can always be achieved by a variant of Algorithm~\ref{alg:gstensor} that uses complex (generalized) Schur decompositions.

Line~\ref{line:gsylvfull} of Algorithm~\ref{alg:gsylvtenrec} is replaced with the following procedure. When $n_{d-1} n_{d} \le n_{\min}^2$,
the matrix
\[ A^\prime_{d-1} = A_d \otimes A_{d-1} 
\]
is formed explicitly. In turn, the $d$-dimensional tensor equation~\eqref{eq:tensorequation3} can equivalently be viewed as the $(d-1)$-dimensional equation
\[
 \bfX^\prime \times_1 A_{1} + \bfX^\prime \times_1 C \times A_{2} \times_3 \cdots  \times_d A^\prime_{d-1} = \bfB^\prime,
\]
with reshaped $\bfX^\prime, \bfB^\prime \in \R^{n_1 \times \cdots \times n_{d-2} \times n_{d-1} n_d}$.
For $d = 2$, this corresponds to the triangular generalized Sylvester equation $A_1 X + C X A_2^T = B$, for which a recursive blocked algorithm has been described in~\cite{Jonsson2002a}.

{\small 
\begin{algorithm}[H] \small
 \caption{Recursive solution of triangular generalized Sylvester equation with tensor structure: $\bfX=\mathtt{mrggsylv} (A_1,\ldots, A_d, C, \bfB)$}
\label{alg:gsylvtenmrg}
\begin{algorithmic}[1]
\IF{$n_1 n_2 \le n_{\min}^2$} \label{line:ifconditiongsylv}
 \STATE Compute $A_{d-1}^\prime = A_{d} \otimes  A_{d-1} $.
 \STATE Reshape $\bfB^\prime(:,\ldots,:, 1:n_{d-1}n_d) = \bfB(:,\ldots,:, 1:n_{d-1},1:n_{d})$.
 \IF{$d = 3$}
 \STATE Solve generalized Sylvester equation $A_1 \bfX^\prime+C\bfX^\prime (A_2^\prime)^T=B^\prime$. \label{line:gsylv}
 \ELSE
 \STATE Call $\bfX^\prime=\mathtt{mrggsylv}(A_1,A_2,\ldots,A_{d-2}, A_{d-1}^\prime,\bfB^\prime)$.
 \ENDIF
 \STATE Reshape $\bfX(:,\ldots,:,1:n_{d-1},1:n_{d}) = \bfX^\prime(:,\ldots,:, 1:n_{d-1} n_{d})$.
\ELSE 
\STATE Determine $n_\mu = \max_{\nu}n_\nu$ and set $k = \lfloor n_\mu/2 \rfloor$. 
\STATE Execute lines \ref{line:startcode} to \ref{line:endcode} of Algorithm~\ref{alg:gsylvtenrec} with calls to {\tt recgsylv} replaced by calls to {\tt mrggsylv}.
\ENDIF
\end{algorithmic}
\end{algorithm}
}
A straightforward extension of the complexity analysis of Algorithm~\ref{alg:laplacemerge} shows that Algorithm~\ref{alg:gsylvtenmrg} requires 
$O\big(n^{d+1} + n_{\min}^{2} n^d)$ flops.

\subsection{Merging dimensions: quasi-triangular case}
\label{subsec:merge_highdim2real}

When using real (generalized) Schur decompositions and, in turn, dealing with upper quasi-triangular coefficients $A_1,\ldots,A_d$, we are facing a situation similar to the one discussed in Section~\ref{sec:quasitriangular}: The merged coefficient matrix $A^\prime_{d-1} = A_d \otimes A_{d-1}$ is, in general, not quasi-triangular. The structure of $A^\prime_{d-1}$ is very similar but not identical with the Laplace-like case. For example, 
comparing~\eqref{eq:a1} with
\begin{equation} \label{eq:a1gsylv}
A_{d-1}={\scriptsize\left[ \arraycolsep=1.4pt\def\arraystretch{1.1}
\begin{array}{cccc}   \times &  \times & \times &  \times \\   0 & \times & \times &  \times \\   0 & \times & \times & \times \\    0 & 0 & 0 & \times \end{array} \right]}, \quad
 A_d={\scriptsize\left[ \arraycolsep=1.4pt\def\arraystretch{1.1} \begin{array}{ccc} \times &  \times & \times \\  0 & \times & \times \\  0 & \times & \times \end{array} \right]},\quad
 A_{d-1}^\prime = {\scriptsize 
\left[ \arraycolsep=1.4pt\def\arraystretch{1.1}
\begin{array}{cccc|cccc|cccc}
\times &  \times & \times &  \times & \times &  \times & \times &  \times & \times &  \times & \times &  \times  \\ 
0 & \times & \times & \times & 0 &  \times & \times & \times & 0 &  \times & \times & \times  \\ 
0 & \times & \times &  \times & 0 & \times & \times &  \times & 0 & \times & \times &  \times  \\ 
0 & 0 & 0 &  \times & 0 & 0 & 0 & \times & 0 & 0 & 0 & \times \\ \hline
0 & 0 & 0 & 0 & \times & \times& \times& \times & \times & \times& \times & \times \\
0 & 0 & 0 & 0 & 0 & \times& \times& \times & 0 & \times& \times& \times \\
0 & 0 & 0 & 0 & 0 & \times & \times & \times& 0 & \times & \times & \times \\
0 & 0 & 0 & 0 & 0 & 0 & 0 & \times& 0& 0 & 0 & \times \\ \hline 
0 & 0 & 0 & 0 & \times & \times& \times& \times & \times & \times& \times& \times  \\
0 & 0 & 0 & 0 & 0 & \times& \times& \times  & 0 & \times& \times& \times  \\
0 & 0 & 0 & 0 & 0 & \times & \times & \times & 0 & \times & \times & \times  \\
0 & 0 & 0 & 0 & 0 & 0 & 0 & \times & 0 & 0 & 0 & \times 
\end{array}
\right]},
\end{equation}
we see that the off-diagonal blocks now have quasi-triangular instead of diagonal structure.
Nevertheless, the properties and techniques discussed in Section~\ref{sec:quasitriangular} carry over verbatim to $A_{d-1}^\prime$. In particular, $A_{d-1}^\prime$ is a block diagonal matrix with diagonal blocks of size at most $2n_{d-1}$. Moreover, a perfect shuffle permutation of the diagonal blocks can again be used to further reduce the size of diagonal blocks. For example, applying this permutation to the second diagonal block of the matrix in~\eqref{eq:a1gsylv} yields:
\[
  P^T A_{d-1}^\prime P = {\scriptsize 
\left[ \arraycolsep=1.4pt\def\arraystretch{1.1}
\begin{array}{cccc|cc|cccc|cc}
\times &  \times & \times &  \times & \times &  \times & \times & \times & \times & \times & \times & \times  \\ 
0 & \times & \times &  \times & 0 & 0 & \times & \times& \times & \times &\times & \times\\ 
0 &  \times & \times &  \times & 0 & 0 & \times & \times & \times & \times & \times & \times  \\ 
0 & 0 & 0 &  \times & 0 & 0 & 0 & 0 & 0 & 0 & \times & \times \\ \hline
0 & 0 & 0 & 0 & \times & \times& \times & \times & \times & \times & \times & \times \\
0 & 0 & 0 & 0 & \times & \times& \times & \times & \times & \times&  \times & \times \\ \hline
0 & 0 & 0 & 0 & 0 & 0 & \times & \times& \times & \times &\times & \times \\
0 & 0 & 0 & 0 & 0 & 0 & \times & \times& \times & \times & \times & \times \\ 
0 & 0 & 0 & 0 & 0 & 0 & \times & \times & \times & \times& \times& \times \\
0 & 0 & 0 & 0 & 0 & 0 & \times & \times & \times & \times& \times & \times  \\ \hline 
0 & 0 & 0 & 0 & 0 & 0 & 0 & 0 & 0 & 0 & \times & \times  \\
0 & 0 & 0 & 0 & 0 & 0 & 0 & 0 & 0 & 0 & \times & \times 
\end{array}
\right].}
 \]This modification  allows to apply Algorithm~\ref{alg:gsylvtenmrg} to quasi-triangular matrices without increased complexity.

\subsection{Numerical experiments}

To give some insight into the performance of Algorithms~\ref{alg:gsylvtenrec} and~\ref{alg:gsylvtenmrg}, we have implemented them in {\sc Matlab} and conducted numerical experiments in the setting described in Section~\ref{sec:numexplaplace}. In particular, we again make use of complex (generalized) Schur decompositions, to avoid that the overhead incurred by the techniques described in Section~\ref{subsec:merge_highdim2real} distorts the picture. To solve the triangular generalized Sylvester equation in Line~\ref{line:gsylv} of Algorithm~\ref{alg:gsylvtenmrg}, we apply {\tt sylvester\_tri} to $A_1 \bfX^\prime E^T+C\bfX^\prime (A_2^\prime)^T=B^\prime$ with $E = I_{n_2}$. 

\begin{figure}
    \centering
    \subfloat[Algorithm~\ref{alg:gsylvtenrec}]{\includegraphics[width=0.49\textwidth]{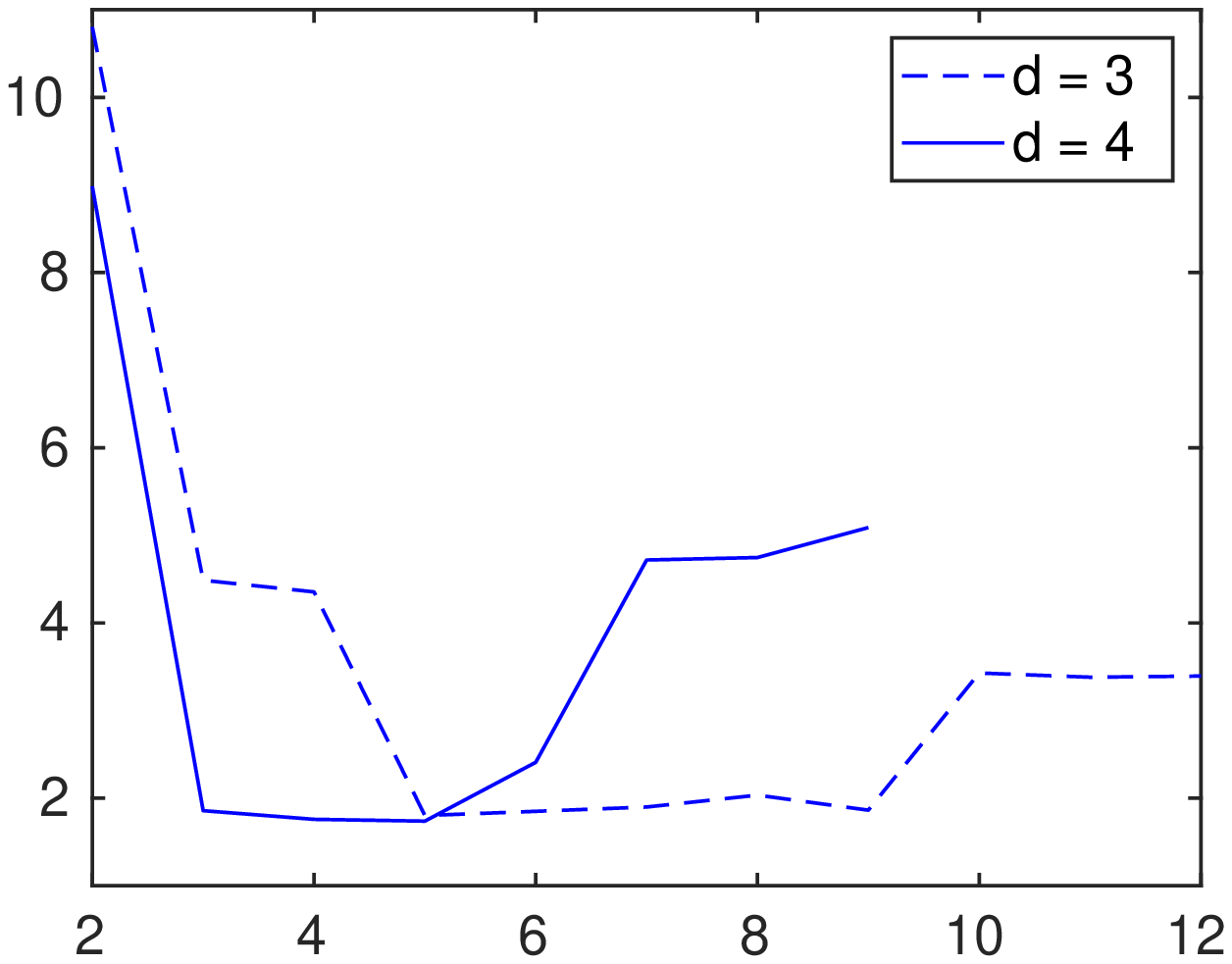}} 
    \subfloat[Algorithm~\ref{alg:gsylvtenmrg}]{\includegraphics[width=0.49\textwidth]{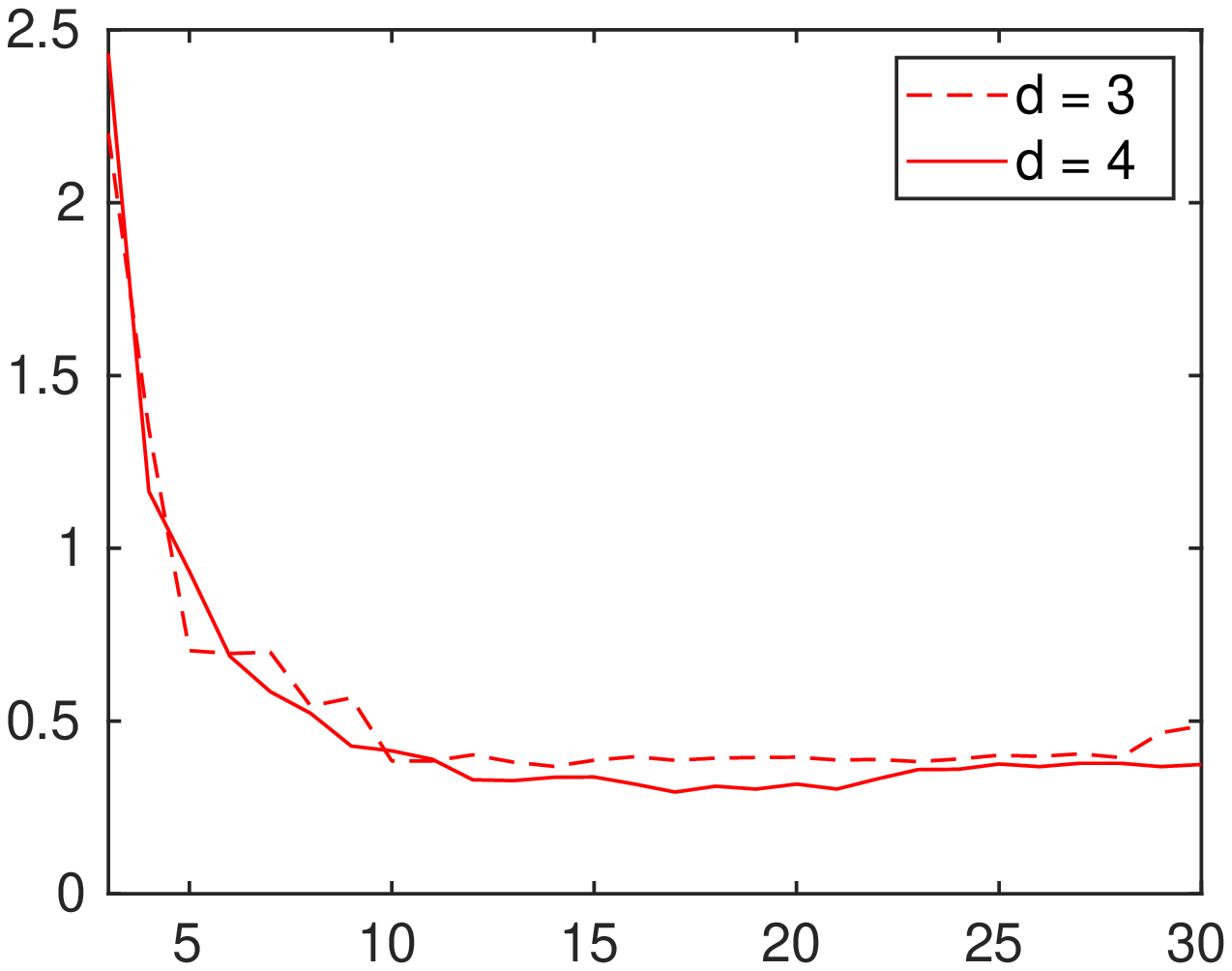}} 
    \caption{Execution times (in seconds) vs. $n_{\min}$ for Algorithms~\ref{alg:gsylvtenrec} and~\ref{alg:gsylvtenmrg} applied to random $n\times\cdots\times n$ tensors with $n = 80$ for $d = 3$ and $n = 25$ for $d = 4$. \label{fig:nmingsylv}}
\end{figure}

\begin{paragraph}{Choice of $n_{\min}$.} Figure~\ref{fig:nmingsylv} shows the performance of Algorithms~\ref{alg:gsylvtenrec} and~\ref{alg:gsylvtenmrg} with respect to the choice of $n_{\min}$. Compared with Algorithms~\ref{alg:laplacerecursive} and~\ref{alg:laplacemerge}, see Figure~\ref{fig:nminlaplace}, the  findings do not differ much. In the following we set $n_{\min} = 8$ for $d = 3$, $n_{\min} = 6$ for $d = 4$ when using Algorithm~\ref{alg:gsylvtenrec}, and $n_{\min} = 15$ for $d = 3$, $n_{\min} = 13$ for $d \ge 4$ when using Algorithm~\ref{alg:gsylvtenmrg}.
\end{paragraph}

\begin{figure}
\begin{minipage}{.5\linewidth}
\centering
\subfloat[$d = 3$]{\includegraphics[scale=.5]{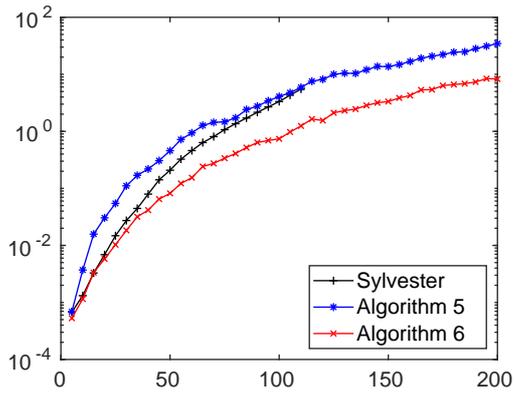}}
\end{minipage}%
\begin{minipage}{.5\linewidth}
\centering
\subfloat[$d = 4$]{\includegraphics[scale=.5]{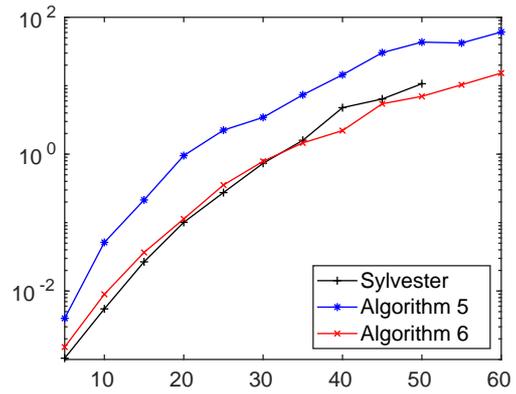}}
\end{minipage} \\
\centering
\subfloat[$d = 5$]{\includegraphics[scale=.5]{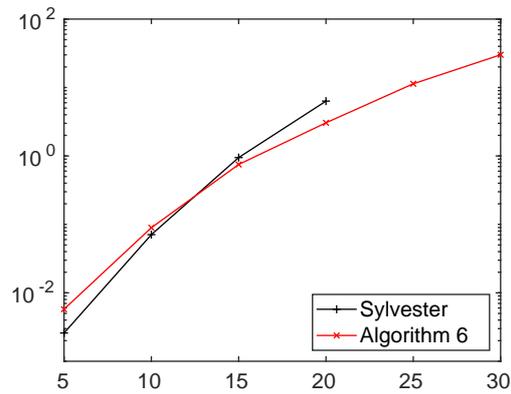}}
    \caption{Execution times (in seconds) vs. $n$ for Algorithms~\ref{alg:gsylvtenrec} and~\ref{alg:gsylvtenmrg} compared to Sylvester solver. \label{fig:compgsylv}}
\end{figure}

\begin{paragraph}{Comparison.}
Figure~\ref{fig:compgsylv} compares the performance of Algorithm~\ref{alg:gsylvtenrec} and Algorithm~\ref{alg:gsylvtenmrg} with the following ``Sylvester solver'':
After reducing the coefficients $A_1,\ldots, A_d,C$ to triangular form and suitably reshaping $\bfB$, one of the Sylvester equations 
\[
\begin{array}{rcl}
 A_1  X + C X (A_3\otimes A_2)^T &=& B \\
 (I \otimes  A_1) X + (A_2 \otimes C ) X (A_4  \otimes A_3 )^T &=& B \\
 (I \otimes  A_1) X + (A_2 \otimes C ) X (A_5\otimes A_4  \otimes A_3 )^T &=& B
 \end{array}
\]
is solved for $d = 3,4,5$ by calling {\tt sylvester\_tri}. 
The results from Figure~\ref{fig:compgsylv} show that Algorithm~\ref{alg:gsylvtenmrg} is always faster than Algorithm~\ref{alg:gsylvtenrec}. The Sylvester solver is slower for sufficiently large $n$; the difference is most pronounced for $d = 3$. Moreover, the Sylvester solver encounters out of memory errors for $n > 110$, $n >50$, $n > 20$ for $d = 3,4,5$, respectively.
\end{paragraph}

